# A multiresolution finite element method based on a new quadrilateral plate element

YiMing Xia

Civil Engineering Department, Nanjing University of

Aeronautics and Astronautics, Nanjing 210016, China

Email:xym4603@sina.com

**ABSTRACT:** A new multiresolution quadrilateral plate element is proposed and a multiresolution finite element method is hence presented. The multiresolution analysis (MRA) framework is formulated out of a mutually nesting displacement subspace sequence, whose basis functions are constructed of scaling and shifting on the element domain of basic node shape function. The basic node shape function is constructed by extending shape function around a specific node. The MRA endows the proposed element with the resolution level (RL) to adjust the element node number, thus modulating structural analysis accuracy accordingly. As a result, the traditional 4-node quadrilateral plate element and method is a monoresolution one and also a special case of the proposed element and method. The meshing for the monoresolution plate element model is based on the empiricism while the RL adjusting for the multiresolution is laid on the rigorous mathematical basis. The accuracy of a structural analysis is actually determined by the RL, not by the mesh. The rational MRA enables the implementation of the multiresolution element method to be more rational and efficient than that of the conventional monoresolution plate element method or other corresponding MRA methods such as the wavelet finite element method, the meshfree method, and the natural element method etc.

## 1. Introduction

Multiresolution analysis (MRA) is a popular technique that has been applied in many domains such as the signal and image processing, the damage detection and health monitoring, the differential equation solution, etc. However, in the field of computational mechanics, the MRA has not been, in a real sense, fully utilized in the numerical solution of engineering problems either by the traditional finite element method (FEM) [1] or by other methods such as the wavelet finite element method (WFEM) [2, 3], the meshfree method (MFM) [4, 5] and the natural element method (NEM) [6, 7] etc.

As is commonly known of the FEM, owing to the invariance of node number a single finite element contains, the finite element can be regarded as a monoresolution one from a MRA point of view and the FEM structural analysis is usually not associated with the MRA concept.

The MRA seems to be rarely used when the FEM is employed to structural analysis. However, it is, in fact, by means of finite element model meshing and re-meshing to modulate analysis accuracy in which a cluster of monoresolution finite elements are assembled together artificially that the rough structural MRA is executed by the FEM. As we can see, in over whole analysis process of a structure by the FEM, there is no solid mathematical foundation for the traditional finite element model meshing and the finite elements are assembled together artificially. The traditional finite element model has to be re-meshed until sufficient accuracy is reached, which leads to the low computation efficiency or convergent rate. The deficiency of the FEM becomes much explicit in the accurate computation of structural problems with local steep gradient such as material nonlinear [8, 9], local damage and crack [10, 11], impacting and exploding problems [12, 13].

The great efforts have been made over the past thirty years to overcome the drawbacks of the FEM with many improved methods to come up, such as WFEM, MFM and NEM etc, which open up a transition from the monoresolution finite element method to the multiresolution finite element method featured with adjustable element node number. Although these MRA methods have illustrated their powerful capability and computational efficiency in dealing with some problems, they always have such major inherent deficiencies as the complexity of shape function construction, the absence of the Kronecker delta property of the shape function and the lack of a solid mathematical basis for the MRA, which make the treatment of element boundary condition complicated and the selection of element node layout empirical, that substantially reduce computational efficiency. Hence, these MRA methods have never found a wide application in engineering practice just as the FEM. In fact, they can be viewed as the intermediate products in the transition of the FEM from the monoresolution to the multiresolution.

The drawbacks of all those MRA methods can be eliminated by the introduction of a new multiresolution finite element method in this paper. With respect to a quadrilateral plate element in the finite element stock, a new multiresolution plate element is formulated by the MRA based on a displacement subspace sequence which is constituted by the translated and scaled version as subspace basis functions of the basic node shape function. The basic node shape function is constructed by extending isoparametrical node shape function for a conventional plate element to other three quadrants around a specific node and the various node shape functions within the multiresolution element can be easily made up of the scaled and shifted version of the basic node shape function. It can be seen that the shape function construction is simple and clear meanwhile the node shape functions hold the Kronecker delta property. In addition, the proposed element method possesses a solid mathematical basis for the MRA, which endows the proposed element with the resolution level (RL) that can be modulated to change the element node number, adjusting structural analysis accuracy accordingly. As a result, the multiresolution quadrilateral plate element method can bring about substantial improvement of the computational efficiency in the structural analysis when compared with the corresponding FEM or other MRA methods.

## 2. The basic node shape function construction

Corresponding to a quadrilateral element of arbitrary shape in a Carstesian coordinate system shown in Fig.1, a isoparametric quadrilateral element in a natural coordinate system shown in Fig.2 is adopted with the transforming relations from the Carstesian coordinate system to the natural one($\xi$ $\eta$) defined as follows:

$$x = \sum_{i=1}^{4} N_i^m x_i \qquad y = \sum_{i=1}^{4} N_i^m y_i \qquad (1)$$

where $x_i$, $y_i$ are the Carstesian coordinate values at the *i*-th node. $N_i^m$ are the conventional shape functions for the different node, which are defined on the domain of $[0,1]^2$ with the bi-linear functions as follows

$$\begin{aligned} N_1^m(\xi,\eta) &= (1-\xi)(1-\eta) \\ N_2^m(\xi,\eta) &= \xi(1-\eta) \\ N_3^m(\xi,\eta) &= \xi\eta \\ N_4^m(\xi,\eta) &= (1-\xi)\eta \end{aligned} \qquad (2)$$

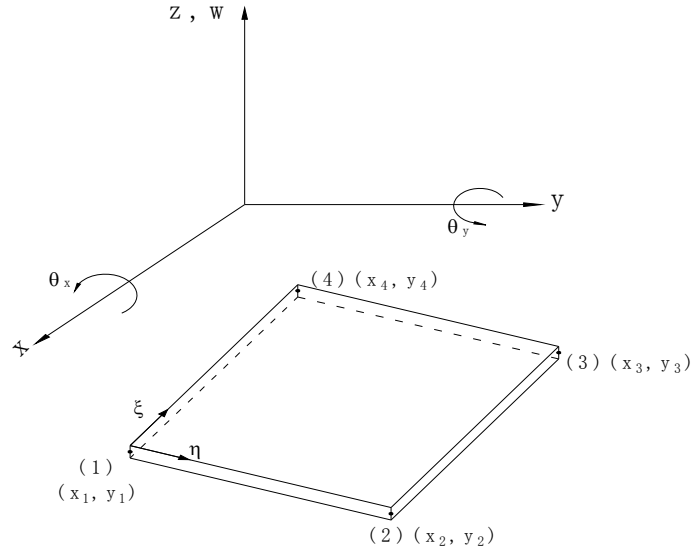

**Fig 1.** A quadrilateral plate bending element

The derivative relationship between the two coordinate systems exists as follows:

$$\left\{ \begin{array}{c} \dfrac{\partial}{\partial x} \\ \dfrac{\partial}{\partial y} \end{array} \right\} = \mathbf{J}_1^{-1} \left\{ \begin{array}{c} \dfrac{\partial}{\partial \xi} \\ \dfrac{\partial}{\partial \eta} \end{array} \right\} \qquad (3)$$

in which $\mathbf{J}_1 = \begin{bmatrix} \partial x/\partial \xi & \partial y/\partial \xi \\ \partial x/\partial \eta & \partial y/\partial \eta \end{bmatrix}$ is the first order Jacobi matrix.

The displacement of a classical bending element shown in Fig.1 can be easily acquired and concisely expressed in terms of natural coordinates as follows [14]:

$$w^e = \sum_{i=1}^{4} N_i^b w_i + \sum_{i=1}^{4} N_{xi}^b \theta_{xi} + \sum_{i=1}^{4} N_{yi}^b \theta_{yi} \qquad (4)$$

where $w^e$ is the transverse displacement in the z axis direction at an arbitrary point of the element. $w_i, \theta_{xi}, \theta_{yi}$ are the transverse, rotational displacements at node $i$ of the element respectively in the Carstesian coordinate system. $N_i^b$, $N_{xi}^b$, $N_{yi}^b$ are the conventional shape functions at the node $i$, which are defined on the domain of $[0,1]^2$ as follows

$$\begin{aligned}
N_1^b(\xi,\eta) &= X_1 Y_1 (X_1 Y_1 - X_2 Y_2 + 2 X_1 X_2 + 2 Y_1 Y_2) \\
N_2^b(\xi,\eta) &= X_2 Y_1 (X_2 Y_1 - X_1 Y_2 + 2 X_1 X_2 + 2 Y_1 Y_2) \\
N_3^b(\xi,\eta) &= X_2 Y_2 (X_2 Y_2 - X_1 Y_1 + 2 X_1 X_2 + 2 Y_1 Y_2) \\
N_4^b(\xi,\eta) &= X_1 Y_2 (X_1 Y_2 - X_2 Y_1 + 2 X_1 X_2 + 2 Y_1 Y_2)
\end{aligned} \qquad (5)$$

$$\begin{aligned}
N_{x1}^b(\xi,\eta) &= X_1 Y_1 (y_{21} X_1 X_2 + y_{41} Y_1 Y_2) \\
N_{x2}^b(\xi,\eta) &= X_2 Y_1 (-y_{21} X_1 X_2 + y_{32} Y_1 Y_2) \\
N_{x3}^b(\xi,\eta) &= X_2 Y_2 (-y_{34} X_1 X_2 - y_{32} Y_1 Y_2) \\
N_{x4}^b(\xi,\eta) &= X_1 Y_2 (y_{34} X_1 X_2 - y_{41} Y_1 Y_2)
\end{aligned} \qquad (6)$$

$$\begin{aligned}
N_{y1}^b(\xi,\eta) &= X_1 Y_1 (-x_{21} X_1 X_2 - x_{41} Y_1 Y_2) \\
N_{y2}^b(\xi,\eta) &= X_2 Y_1 (x_{21} X_1 X_2 - x_{32} Y_1 Y_2) \\
N_{y3}^b(\xi,\eta) &= X_2 Y_2 (x_{34} X_1 X_2 + x_{32} Y_1 Y_2) \\
N_{y4}^b(\xi,\eta) &= X_1 Y_2 (-x_{34} X_1 X_2 + x_{41} Y_1 Y_2)
\end{aligned} \qquad (7)$$

where $X_1 = 1-\xi, X_2 = \xi, Y_1 = 1-\eta, Y_2 = \eta, x_{ij} = x_i - x_j, y_{ij} = y_i - y_j$ ($i$=1~4, $j$=1~4).

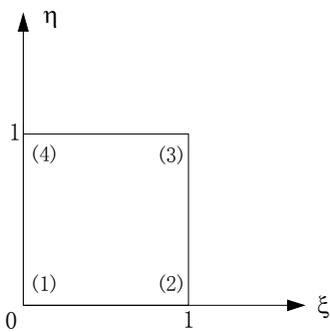

**Fig 2.** The node shape function domain of a classical quadrilateral plate element

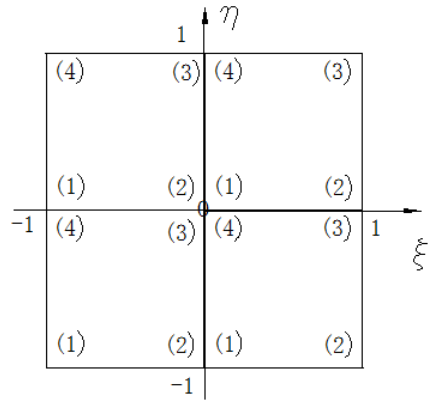

**Fig 3.** The extended shape function domain for the node (1)

As to the proposed element, the shape functions, regarding the node (1) as shown in Fig 2. initially defined on the domain of [0, 1]$^2$, should be extended to the domain of $[-1,1]^2$ by means of shifting the element around the node (1) vertically, horizontally and obliquely respectively to the other three quadrants, thus covering the eight nodes adjacent to the node (1) as shown in Fig 3, the basic shape function for the node (1) at the point coordinate of (0, 0) can finally be defined as follows:

$$\phi_1(\xi,\eta) := \begin{cases} N_1^b(\xi,\eta) & \xi \in [0,1], \eta \in [0,1] \\ N_2^b(1+\xi,\eta) & \xi \in [-1,0], \eta \in [0,1] \\ N_3^b(1+\xi,1+\eta) & \xi \in [-1,0], \eta \in [-1,0] \\ N_4^b(\xi,1+\eta) & \xi \in [0,1], \eta \in [-1,0] \\ 0 & \xi \notin [-1,1], \eta \notin [-1,1] \end{cases} \quad (8)$$

$$\phi_2(\xi,\eta) := \begin{cases} N_{x1}^b(\xi,\eta) & \xi \in [0,1], \eta \in [0,1] \\ N_{x2}^b(1+\xi,\eta) & \xi \in [-1,0], \eta \in [0,1] \\ N_{x3}^b(1+\xi,1+\eta) & \xi \in [-1,0], \eta \in [-1,0] \\ N_{x4}^b(\xi,1+\eta) & \xi \in [0,1], \eta \in [-1,0] \\ 0 & \xi \notin [-1,1], \eta \notin [-1,1] \end{cases} \quad (9)$$

$$\phi_3(\xi,\eta) := \begin{cases} N_{y1}^b(\xi,\eta) & \xi \in [0,1], \eta \in [0,1] \\ N_{y2}^b(1+\xi,\eta) & \xi \in [-1,0], \eta \in [0,1] \\ N_{y3}^b(1+\xi,1+\eta) & \xi \in [-1,0], \eta \in [-1,0] \\ N_{y4}^b(\xi,1+\eta) & \xi \in [0,1], \eta \in [0,-1] \\ 0 & \xi \notin [-1,1], \eta \notin [-1,1] \end{cases} \quad (10)$$

With help of eq.(3), it can be seen that the Kronecker delta property holds for the basic node shape functions $\phi_1(\xi,\eta)$, $\phi_2(\xi,\eta)$, $\phi_3(\xi,\eta)$ with respect to the Carstesian coordinate system:

$$\begin{cases} \phi_1(0,0)=1 & \phi_1(Anodes)=0 & \dfrac{\partial \phi_1(0,0)}{\partial x}=0 \\ \dfrac{\partial \phi_1(0,0)}{\partial y}=0 & \dfrac{\partial \phi_1(Anodes)}{\partial x}=0 & \dfrac{\partial \phi_1(Anodes)}{\partial y}=0 \\ \phi_2(0,0)=0 & \phi_2(Anodes)=0 & \dfrac{\partial \phi_2(0,0)}{\partial x}=0 \\ \dfrac{\partial \phi_2(0,0)}{\partial y}=1 & \dfrac{\partial \phi_2(Anodes)}{\partial x}=0 & \dfrac{\partial \phi_2(Anodes)}{\partial y}=0 \\ \phi_3(0,0)=0 & \phi_3(Anodes)=0 & \dfrac{\partial \phi_3(0,0)}{\partial x}=-1 \\ \dfrac{\partial \phi_3(0,0)}{\partial y}=0 & \dfrac{\partial \phi_3(Anodes)}{\partial x}=0 & \dfrac{\partial \phi_3(Anodes)}{\partial y}=0 \end{cases} \quad (11)$$

where the notation Anodes is referred to the eight nodes adjacent to the basic node (1)

The basic node shape functions $\phi_1(\xi,\eta), \phi_2(\xi,\eta), \phi_3(\xi,\eta)$ in the natural coordinate system are illustrated in Figs.4 *a,b,c* respectively

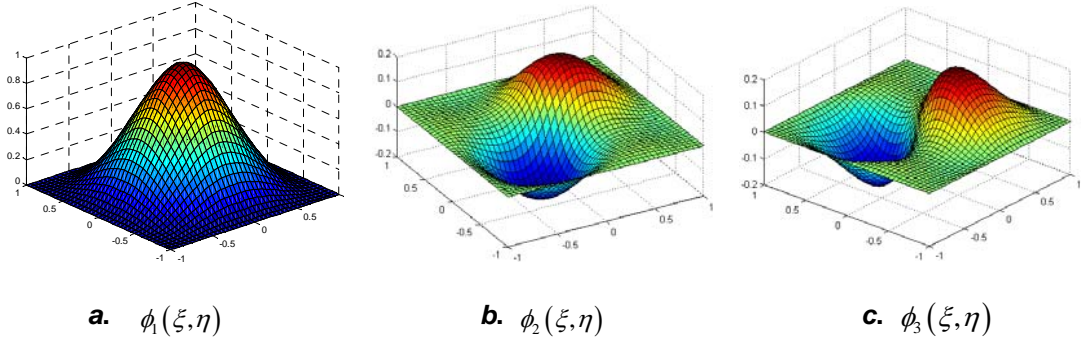

**a.** $\phi_1(\xi,\eta)$    **b.** $\phi_2(\xi,\eta)$    **c.** $\phi_3(\xi,\eta)$

**Fig 4.** The basic node shape functions $\phi_1(\xi,\eta), \phi_2(\xi,\eta), \phi_3(\xi,\eta)$ on the domain [-1,1]× [-1,1]

Obviously, the basic shape functions $\phi_1(\xi,\eta), \phi_2(\xi,\eta), \phi_3(\xi,\eta)$ are continuous respectively.

## 2  Basis function construction for a displacement subspace sequence forming a new MRA

In order to carry out a MRA of a thin plate structure, the mutual nesting displacement subspace sequence for a plate element should be established. In this paper, a totally new technique is proposed to construct the MRA which is based on the concept that a subspace sequence (multi-resolution subspaces) can be formulated by subspace basis function vectors at different resolution levels whose elements-scaling function vector can be constructed by scaling and shifting on the domain $[0,1]^2$ of the basic node shape functions. As a result, the displacement subspace basis function vector at an arbitrary resolution level (RL) of $(m+1)\times(n+1)$ for a quadrilateral plate element is formulated as follows:

$$\mathbf{\Psi}_{mn} = \begin{bmatrix} \mathbf{\Phi}_{mn,00} & \cdots & \mathbf{\Phi}_{mn,rs} & \cdots & \mathbf{\Phi}_{mn,mn} \end{bmatrix} \tag{12}$$

where $\mathbf{\Phi}_{mn,rs} = \begin{bmatrix} \phi_{1mn,rs} & \phi_{2mn,rs}/n & \phi_{3mn,rs}/m \end{bmatrix}$ is the scaling basis function vector, $\phi_{1mn,rs} = \phi_1(m\xi - r, n\eta - s)$, $\phi_{2mn,rs} = \phi_2(m\xi - r, n\eta - s)$, $\phi_{3mn,rs} = \phi_3(m\xi - r, n\eta - s)$, $m$, $n$ denoted as the positive integers, the scaling parameters in $\xi, \eta$ directions respectively. $r$, $s$ as the positive integers, the node position parameters, that is $r = 0,1,2,3\cdots m$, $s = 0,1,2,3\cdots n$, Here, $(m\xi - r) \in [-1,1]$ $(n\eta - s) \in [-1,1]$, $\xi \in [0,1], \eta \in [0,1]$.

It is seen from Eq. (12) that the nodes for the scaling process are equally spaced on the domain $[0,1]^2$ with a step size of $1/m$ in $\xi$ and $1/n$ in $\eta$ directions respectively.

Scaling of the basic node shape functions on the domain of $[-1, 1]^2$ (precisely on the domain of $\left[-\dfrac{1}{m}, \dfrac{1}{m}\right] \times \left[-\dfrac{1}{n}, \dfrac{1}{n}\right]$) and then shifting to other nodes $\left(\dfrac{r}{m}, \dfrac{s}{n}\right)$ on the element domain of $[0, 1]^2$ will produce the various node shape functions that are shown in the Fig 5. at the RL of $2\times 2$, $3\times 3$.

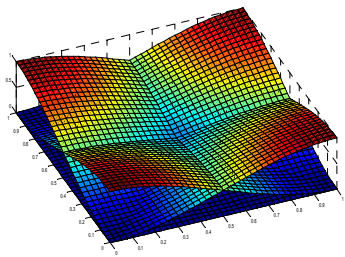 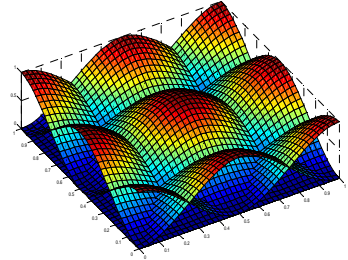

**a.** $\phi_1(\xi,\eta)$ （RL=2×2）　　　　　　**b.** $\phi_1(\xi,\eta)$ （RL=3×3）

**Fig 5**. The scaled and shifted version of the basic node shape function $\phi_1(\xi,\eta)$ on the domain

$$[0,1]\times[0,1]$$

Since the elements in the basis functions are linearly independent with the various scaling and the different shifting parameters, the subspaces in the subspace sequence can be established and are mutually nested, thus formulating a MRA framework, that is

$$\begin{cases} \mathbf{W}_{mn} = \begin{bmatrix} V_{11} \cdots & V_{ij} \cdots & V_{mn} \end{bmatrix} \\ V_{ij} := \mathrm{span}\{\mathbf{\Psi}_{ij} : i, j \in Z\} \\ V_{ij} \subset V_{i(j+1)} \quad V_{ij} \subset V_{(i+1)j} \quad V_{ij} \subset V_{(i+1)(j+1)} \end{cases} \tag{13}$$

where $Z$ denoted as the positive integers, $V_{ij}$ as displacement subspace at the resolution level of $(i+1)\times(j+1)$.

Thus, it can be found that the mutually nesting displacement subspace sequence $\mathbf{W}_{mn}$ can be taken for a solid mathematical foundation for the MRA framework and $V_{11}$ is equivalent to the displacement field for a traditional 4-node plate element that is the reason why the traditional quadrilateral plate element is regarded as a mono-resolution one and also a special case of the multiresolution quadrilateral plate element.

Based the MRA established, the deflection of a quadrilateral plate element in the displacement subspace at RL of $(m+1)\times(n+1)$ can be defined as follows

$$w_{mn}^e = \mathbf{\Psi}_{mn} \mathbf{a}_{mn}^e \tag{14}$$

where $\mathbf{a}_{mn}^e = \left[\left[w_{00}, \theta_{x00}, \theta_{y00}\right] \ldots \left[w_{rs}, \theta_{xrs}, \theta_{yrs}\right] \ldots \left[w_{mn}, \theta_{xmn}, \theta_{ymn}\right]\right]^T$, $w_{rs}, \theta_{xrs}, \theta_{yrs}$ are the transverse and rotational displacements respectively at the element node $\left(\dfrac{r}{m}, \dfrac{s}{n}\right)$.

It is obvious that the proposed multi-resolution element is a meshfree one whose nodes are uniformly scattered, node number and position fully determined by the *RL*. When the scaling parameter $m=n=1$ ($RL=2\times2$), that is a traditional 4-node quadrilateral plate element, eq. (14) will be reduced to eq. (4).

## 3  Multiresolution quadrilateral plate element formulation

According to the classical assumption of the thin plate theory, the generalized function of potential energy in a displacement subspace at the resolution level $(m+1)\times(n+1)$ for a quadrilateral plate element is

$$\Pi(V_{mn}) = \frac{1}{2}\int_0^a\int_0^b [\kappa]_{mn}^T [D_b][\kappa]_{mn}\,dxdy - \int_0^a\int_0^b q w_{mn}^e\,dxdy - \sum_i Q_i w_{mni}^e \tag{15}$$

where $[\kappa]_{mn} = -\begin{bmatrix} \dfrac{\partial^2 w_{mn}^e}{\partial x^2} \\ \dfrac{\partial^2 w_{mn}^e}{\partial y^2} \\ 2\dfrac{\partial^2 w_{mn}^e}{\partial x \partial y} \end{bmatrix}$, $[D_b] = C_b \begin{bmatrix} 1 & \mu & 0 \\ \mu & 1 & 0 \\ 0 & 0 & (1-\mu)/2 \end{bmatrix}$, $C_b = \dfrac{Eh^3}{12(1-\mu^2)}$, $E$ is

the material Young modulus, $h$ the thickness of the element, $\mu$ the Poisson's ratio, $q$ distributed transverse loadings, $Q$ the lump transverse loadings.

$$[\kappa]_{mn} = [\mathbf{B}_{00}, \cdots \mathbf{B}_{rs}, \cdots \mathbf{B}_{mn}] \mathbf{a}_{mn}^e \tag{16}$$

where

$$\mathbf{B}_{rs} = -\left[\frac{\partial^2 \mathbf{\Phi}_{mn,rs}}{\partial x^2} \quad \frac{\partial^2 \mathbf{\Phi}_{mn,rs}}{\partial y^2} \quad 2\frac{\partial^2 \mathbf{\Phi}_{mn,rs}}{\partial x \partial y}\right]^T$$

$$\begin{Bmatrix} \partial^2 \mathbf{\Phi}_{mn,rs}/\partial x^2 \\ \partial^2 \mathbf{\Phi}_{mn,rs}/\partial y^2 \\ \partial^2 \mathbf{\Phi}_{mn,rs}/\partial x \partial y \end{Bmatrix} = \mathbf{J}_2^{-1} \begin{Bmatrix} \partial^2 \mathbf{\Phi}_{mn,rs}/\partial \xi^2 \\ \partial^2 \mathbf{\Phi}_{mn,rs}/\partial \eta^2 \\ \partial^2 \mathbf{\Phi}_{mn,rs}/\partial \xi \partial \eta - \alpha' \partial \mathbf{\Phi}_{mn,rs}/\partial \xi + \beta' \partial \mathbf{\Phi}_{mn,rs}/\partial \eta \end{Bmatrix},$$

$$\mathbf{J}_2 = \begin{bmatrix} (\partial x/\partial \xi)^2 & (\partial y/\partial \xi)^2 & 2(\partial x/\partial \xi)(\partial y/\partial \xi) \\ (\partial x/\partial \eta)^2 & (\partial y/\partial \eta)^2 & 2(\partial x/\partial \eta)(\partial y/\partial \eta) \\ (\partial x/\partial \xi)(\partial x/\partial \eta) & (\partial y/\partial \xi)(\partial y/\partial \eta) & (\partial x/\partial \xi)(\partial y/\partial \eta)+(\partial y/\partial \xi)(\partial x/\partial \eta) \end{bmatrix}$$

is the second order Jacobi matrix. $\alpha' = (\alpha \cdot \partial y/\partial \eta - \beta \cdot \partial x/\partial \eta)/|\mathbf{J}_1|$,

$\beta' = (\alpha \cdot \partial y/\partial \xi - \beta \cdot \partial x/\partial \xi)/|\mathbf{J}_1|$, $\alpha = x_1 - x_2 + x_3 - x_4$ $\beta = y_1 - y_2 + y_3 - y_4$.

Substitute Eq.(14), Eq.(16) into Eq.(15), the concise expression can be obtained after reassembling as follows:

$$\Pi_p(V_{mn}) = \frac{1}{2}\mathbf{a}_{mn}^{eT}\mathbf{K}_{mn}^{e}\mathbf{a}_{mn}^{e} - \mathbf{a}_{mn}^{eT}\mathbf{f}_{mn}^{e} - \mathbf{a}_{mn}^{eT}\mathbf{F}_{mn}^{e} \qquad (17)$$

where $\mathbf{K}_{mn}^{e}$ is the element stiffness matrix, $\mathbf{f}_{mn}^{e}$ the element distributed loading column vector, $\mathbf{F}_{mn}^{e}$ the element lump loading column vector.

According to the potential energy minimization principle, let $\delta\Pi_p(V_{mn}) = 0$, the plate element equilibrium equations can be obtained as follows,

$$\mathbf{K}_{mn}^{e}\mathbf{a}_{mn}^{e} = \mathbf{f}_{mn}^{e} + \mathbf{F}_{mn}^{e} \qquad (18)$$

The element expression of the stiffness matrix $\mathbf{K}_{mn}^{e}$, and the loading column vectors $\mathbf{f}_{mn}^{e}$, $\mathbf{F}_{mn}^{e}$ can be given as follows:

$$\mathbf{K}_{mn}^{e} = \begin{bmatrix} \mathbf{k}_{00}^{00} & \cdots & \mathbf{k}_{rs}^{00} & \cdots & \mathbf{k}_{mn}^{00} \\ \cdot & & \cdot & & \cdot \\ \cdot & & \cdot & & \cdot \\ \cdot & & \cdot & & \cdot \\ \mathbf{k}_{00}^{rs} & \cdots & \mathbf{k}_{rs}^{rs} & \cdots & \mathbf{k}_{ij}^{rs} \\ \cdot & & \cdot & & \cdot \\ \cdot & & \cdot & & \cdot \\ \cdot & & \cdot & & \cdot \\ \mathbf{k}_{00}^{mn} & \cdot & \mathbf{k}_{rs}^{mn} & \cdot & \mathbf{k}_{mn}^{mn} \end{bmatrix} \qquad (19)$$

where the superscript denoted as the row number of the matrix and the subscript as the aligned element node numbering ($r$, $s$). In terms of the properties of the node shape functions, we have

$$\begin{cases} \mathbf{k}_{rs}^{rs} = \sum_{\substack{|c-r|\leq 1 \\ |d-s|\leq 1}} \mathbf{k}_{cd,rs} \\ \mathbf{k}_{rs}^{rs} = \mathbf{k}_{cd,rs} = 0, when |c-r|>1, |d-s|>1 \end{cases} \quad (20)$$

in which $\mathbf{k}_{cd,rs}$ is the coupled node stiffness matrix relating the node ($c$, $d$) to ($r$, $s$).

$$\mathbf{k}_{cd,rs} = \int_0^1 \int_0^1 \mathbf{B}_{cd}^T [D_b] \mathbf{B}_{rs} |\mathbf{J}_1| d\xi d\eta \quad (21)$$

$$\begin{cases} \mathbf{f}_{mn}^e = \int_0^1 \int_0^1 \mathbf{\Psi}_{mn}^T q |\mathbf{J}_1| d\xi d\eta \\ \mathbf{F}_{mn}^e = \mathbf{\Psi}_{mn}^T \mathbf{P} \end{cases} \quad (22)$$

where $\mathbf{\Psi}_{mn}$ is the shape function matrix, $\mathbf{P}$ is the lump loading vector.

## 4 Transformation matrix

In order to carry out structural analysis, the element stiffness $\mathbf{K}_{mn}^e$, the loading column vectors $\mathbf{f}_{mn}^e$, $\mathbf{F}_{mn}^e$ should be transformed from the element local coordinate system (*xyz*) to the structural global coordinate system (*XYZ*). The transforming relations from the local to the global are defined as follows:

$$\mathbf{K}_{mnl}^i = \mathbf{T}_{mn}^{eT} \mathbf{K}_{mn}^e \mathbf{T}_{mn}^e \quad (23)$$

$$\mathbf{f}_{mn}^i = \mathbf{T}_{mnl}^{eT} \mathbf{f}_{mn}^e \quad (24)$$

$$\mathbf{F}_{mn}^i = \mathbf{T}_{mnl}^{eT} \mathbf{F}_{mn}^e \quad (25)$$

where $\mathbf{K}_{mn}^i$ is the element stiffness matrix, $\mathbf{f}_{mn}^i$, $\mathbf{F}_{mn}^i$ the element loading column vectors under the global coordinate system. $\mathbf{T}_{mn}^e$ is the element transformation matrix defined as follows:

$$\mathbf{T}_{mn}^e = \begin{bmatrix} \lambda_{11} & & & & \mathbf{0} \\ & \ldots & & & \\ & & \lambda_{ij} & & \\ & & & \ldots & \\ \mathbf{0} & & & & \lambda_{mn} \end{bmatrix} \quad \lambda_{ij} = \begin{bmatrix} \cos\theta_{zZ} & 0 & 0 \\ 0 & \cos\theta_{xX} & \cos\theta_{xY} \\ 0 & \cos\theta_{yX} & \cos\theta_{yY} \end{bmatrix}$$

where $\theta$ is the intersection angle between the local and the global coordinate axes.

The structural global stiffness $\mathbf{K}_{mn}$, and the global loading column vectors $\mathbf{f}_{mn}$, $\mathbf{F}_{mn}$ can be obtained by splicing $\mathbf{K}^i_{mn}$ $\mathbf{f}^i_{mn}$, $\mathbf{F}^i_{mn}$ of the element respectively

## 5  Numerical example

**Example 1.** As shown in Fig.6, a two opposite edge simply supported and other two free $60^0$ skew plate with the geometric configuration of length *L* and the Poisson's ratio $\mu = 0.3$ is subjected to the uniform transverse loading of magnitude *q*. Evaluate the deflection at the center point of the plate.

The displacement responses are found by the proposed quadrilateral plate element model, the traditional 4-node quadrilateral plate element model and the wavelet element model based on two-dimensional tensor product B-spline wavelet on the interval (BSWI) [3] respectively. The BSWI is chosen because it is the best one among all existing wavelets in approximation of numerical calculation [15] and directly constructed by the tensor product of the wavelets expansions at each coordinate. The central deflections of the plate are summarized in table.1.

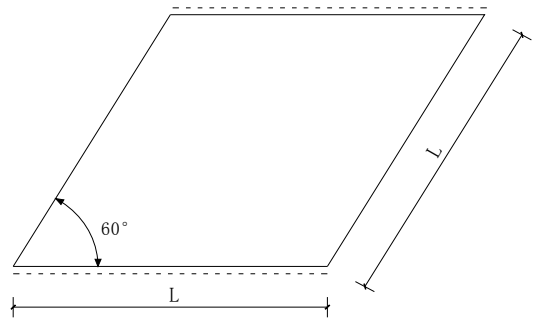

**Fig. 6.**  A skew plate

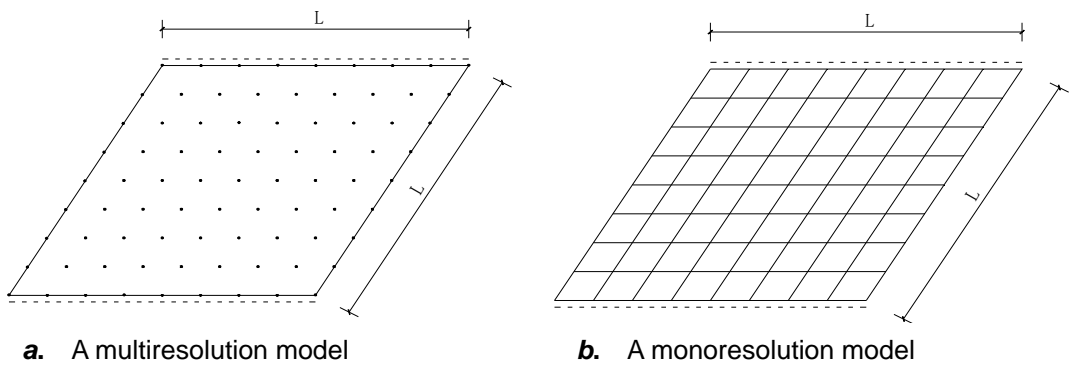

*a.*  A multiresolution model    *b.*  A monoresolution model

**Fig. 7**  The finite element model for the skew plate

Table.1. the center point deflection ($w/\ qL^4/100D_0$)

| Element type | | deflection |
|---|---|---|
| The proposed (RL) | The conventional (mesh) | |
| 9×9 | 8×8 | 0.7876 |
| 13×13 | 12×12 | 0.7909 |
| 17×17 | 16×16 | 0.7920 |
| One BSWI [3] | | 0.7925 |
| Analytical [16] | | 0.7945 |

   The multiresolution models composed of one proposed multiresolution quadrilateral plate element with the RLs of 9×9, 13×13 and 17×17, as shown in Fig.7.*a* are adopted, the monoresolution models composed of meshes of 8×8, 12×12 and 16×16, as displayed in Fig.7.*b*, are also employed and the wavelet models made up of one 2D BSWI element (B-spline wavelet on the interval element) of the jth scale=3, the mth order =4 are used respectively abbreviated as BSWI43 with the DOF of 11×11. The RL of each proposed and the corresponding meshes of the conventional are compared. It can be seen that the analysis accuracies with the proposed quadrilateral element are gradually improved respectively with the RL reaching high. Although the BSWI43 is of high accuracy, when compared with the proposed, the deficiencies of the BSWI element are obvious as follows. In light of tensor product formulation of the multidimensional MRA framework [2,15], the DOF of a multi-dimensional BSWI element will be so drastically increased from that of a one-dimensional element in an irrational way, resulting in complex shape functions and substantial reduction of the computational efficiency. Secondly, due to the absence of Kronecker delta property of the tensor-product constructed shape functions, the special treatments should be taken to deal with the element boundary condition, which will bring about low computational efficiency. Thirdly, there exists no such a parameter as the RL with a clear mathematical sense. In addition, the RLs of the proposed and the corresponding meshes of the conventional are displayed in Table1. It can be found that the analysis accuracies with the proposed and the conventional are gradually improved respectively with the RL reaching high and the mesh approaching dense. However, the RL adjusting is more rationally and efficiently to be implemented than the meshing and the re-meshing for the following two reasons. Firstly, the RL adjusting is based on the MRA framework that is constructed on a solid mathematical basis while the meshing or remeshing, which resorts to the empiricism, has no MRA framework. Secondly, the stiffness matrix and the loading column vectors of the proposed element can be obtained automatically around the nodes while those of the traditional 4-node quadrilateral plate elements obtained by the artificially complex reassembling around the elements. Thus, the computational efficiency of the proposed element method is higher than the traditional one. In this way, the proposed plate element exhibits its strong capability of accuracy adjustment and its high power of resolution to identify details (nodes) of deformed structure by means of modulating its resolution level, just as a multiresolution camera with a pixel in its taken photo as a node in the proposed element. There appears no mesh in the proposed element just as no grid in the photo. Thus, an element of superior analysis accuracy surely has more nodes when compared with that of the inferior just as a clearer photo contains more pixels.

**Example 2.** A circular ring slab is subjected to the uniform transverse loading *q* as shown in Fig.8 with its boundary conditions as: the inner edge is free and the outer edge is fully clamped, and its geometry and physical parameters as: the inner radius *b*, the outer radius *a*, the thickness *t*, the elasticity modulus *E*, the Poisson's ratio $\mu$. Find the displacement around the inner free-edge ring of the slab.

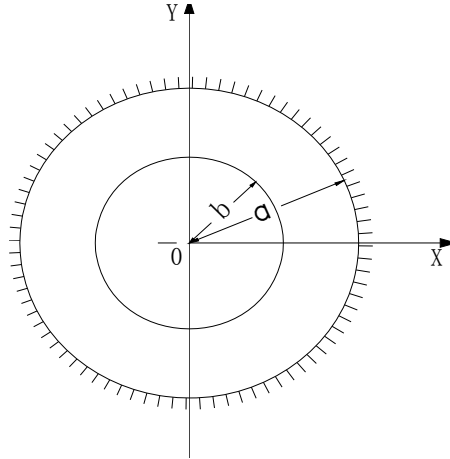

**Fig. 8.** A circular ring slab

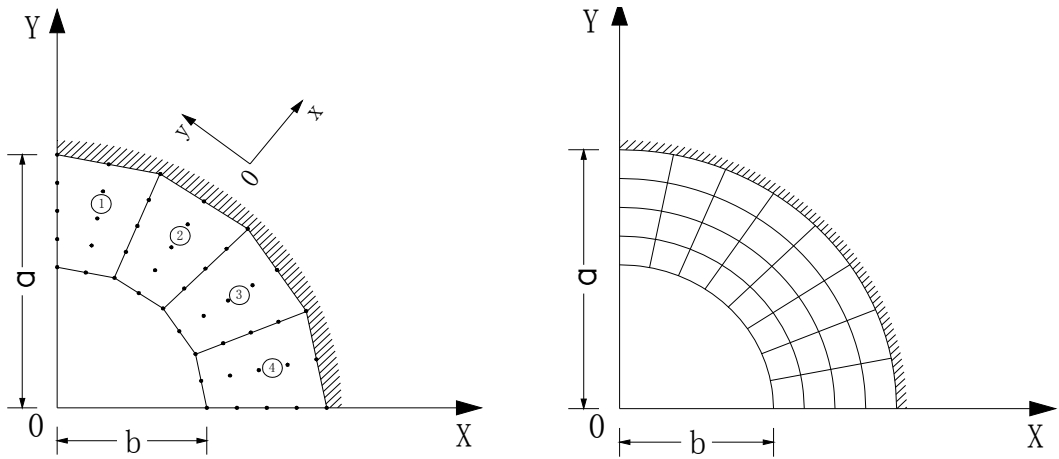

*a* A.multiresolution model        *b.* A monoresolution model

**Fig. 9.** The finite element model for the 1/4 circular ring slab

To calculate the displacement responses, symmetry conditions are exploited and only the first quadrant of the plate is discretized. the mutiresolution plate model is herein composed of four multiresolution quadrilateral plate elements ①,②,③,④ with the element RL of 5×3, hence the 1/4 slab RL of 5×9, as shown in Fig. 9*a*., and the monoresolution one composed of the mesh of 4×8 displayed in Fig.9*b*. In the analysis process, these four multiresolution elements are spliced together along the common intersection boundary and the analysis accuracy can be modulated by means of adjusting the RL. With respect to the conventional monoresolution, the structure is meshed into a group of monoresolution elements and the analysis accuracy is improved only by means of re-meshing. It can be seen that the RL adjusting is more rationally and easily to be implemented than the re-meshing because the proposed multiresolution element model of the circular ring plate structure contains much less elements than the monoresolution, hence requiring much less times of the transformation matrix multiplying, which results in much higher computational efficiency for the proposed element method than that for the traditional element method. The displacement around the inner free-edge ring of the slab is summarized in Table.2. The RL of the intersection boundary should be the same as that of the adjacent element just as PS (Photoshop) four photos.

**Table.2**. The maximum displacement ($w/qa^4/Et^3$) around the inner free-edge ring of the slab

| Element type | | deflection |
|---|---|---|
| The proposed (RL) | The conventional (mesh) | |
| 5×9 | 4×8 | 0.0576 |
| Analytical [17] | | 0.0575 |

## 6  Discussion

From the two numerical examples above, it is shown that based on the multiresolution quadrilateral plate element formulation, a new multiresolution finite element method is introduced, which incorporates such main steps as RL adjusting, element matrix formation, element matrix transformation from a local coordinate system to a global one and global structural matrix formation by splicing of the element matrices. Owing to the existence of the new MRA framework, the RL adjusting for the proposed method is more rationally and easily to be implemented than the meshing and re-meshing for the traditional 4-node quadrilateral plate element method. Due to the basic node shape function, the stiffness matrix and the loading column vectors of a proposed element can be automatically acquired through quadraturing around nodes in the element matrix formation step while those of the traditional 4-node quadrilateral plate obtained through complex artificially reassembling of the element matrix around the elements in the re-meshing process, which contributes a lot to computation efficiency improvement of the proposed method. Moreover, since the multiresolution quadrilateral plate element model of a structure usually contains much less elements than the traditional monoresolution element model, thus requiring much less times of transformation matrix multiplying, the computation efficiency of the proposed method appears much higher than the traditional in the step of element matrix transformation. In addition, because of the simplicity and clarity of the shape function formulation with the Kronecker delta property, and the solid mathematical basis of the new MRA framework, the proposed method is also superior to other corresponding MRA methods in terms of the computational efficiency, the application flexibility and extent. Hence, taking all those causes into account, the conclusion can be drawn that the multiresolution quadrilateral plate element method is more rationally, easily and efficiently to be executed, when compared with the traditional 4-node quadrilateral plate element method or other corresponding MRA methods, and the proposed plate element is the most accurate one formulated ever since.

## 7  Conclusions

A new multiresolution finite element method that has both high power of resolution and strong flexibility of analysis accuracy is introduced into the field of numerical analysis. The method possesses such prominent features as follows:

1. A novel technique is proposed to construct the simple and clear basic node shape function that holds Kronecker delta property.
2. A mathematical basis for the MRA framework, that is the mutually nesting displacement subspace sequence, is constituted out of the scaled and shifted version of the basic node shape function. The MRA framework endows the plate element with the RL to adjust the element node number, modulating the analysis accuracy of structure accordingly. Hence, the traditional 4-node quadrilateral plate element and method is a monoresolution one and also a special case of the proposed. An element of superior analysis accuracy surely contains more nodes when compared with that of the inferior.
3. The RL adjusting for the multiresolution plate element model is laid on the rigorous mathematical basis while the meshing or remeshing for the monoresolution is based on the empiricism. Hence, the implementation of the proposed element method is more rational and efficient than that of the traditional or other MRA methods such as the wavelet finite element method, the meshfree method, and the natural element method etc.
4. A quite new concept is introduced into the FEM that the structural analysis accuracy is actually determined by the RL-the density of node uniform distribution, not by the mesh.
5. With advent of the multiresolution finite element method, the rational MRA will find a wide application in numerical solution of engineering problems in a real sense.

The upcoming work will be focused on the treatment of interface between multiresolution elements of different RL. The interface may be extended to the bridging domain in which the transitional element (expanded Serendipity element) could be used just as PS photos of different RL.

## 8 Acknowledgement

The author would like to thank the referees for their valuable comments, also Prof ShaoLin Chen and associate Prof Gan Tan for their assistance.

## Appendix 1

The equivalent node loading $\mathbf{f}$ of a uniform distributed loading $q_0$ over, or the equivalent node loading $\mathbf{F}$ of a lump loading P applied at the centre of a cell, which is a 4-node subdomain in an element, can be read as follows:

$$\mathbf{F} = P \begin{bmatrix} 1/4 & (y_{21}+y_{41})/16 & (-x_{21}-x_{41})/16 & 1/4 & (-y_{21}+y_{32})/16 & (x_{21}-x_{32})/16 \\ 1/4 & (y_{34}-y_{41})/16 & (-x_{34}+x_{41})/16 & 1/4 & (-y_{34}-y_{32})/16 & (x_{34}+x_{32})/16 \end{bmatrix}^T$$

$$\mathbf{f} = (q_0/180) \begin{bmatrix} \overline{Z}_1 & \overline{T}_{x1} & \overline{T}_{y1} & \overline{Z}_2 & \overline{T}_{x2} & \overline{T}_{y2} \\ \overline{Z}_3 & \overline{T}_{x3} & \overline{T}_{y3} & \overline{Z}_4 & \overline{T}_{x4} & \overline{T}_{y4} \end{bmatrix}^T$$

$\overline{Z}_1 = 36A - 36B + 180C$, $\overline{T}_{x1} = y_{21}(3A - 5B - 30C) + y_{41}(5A - 3B - 30C)$, $\overline{T}_{y1} = -x_{21}(3A - 5B - 30C) - x_{41}(5A - 3B - 30C)$, $\overline{Z}_2 = -36A - 36B + 180C$, $\overline{T}_{x2} = y_{21}(3A + 5B - 30C) - y_{32}(5A + 3B - 30C)$, $\overline{T}_{y2} = -x_{21}(3A + 5B - 30C) - x_{32}(5A + 3B - 30C)$, $\overline{Z}_3 = -36A + 36B + 180C$, $\overline{T}_{x3} = -y_{34}(3A - 5B - 30C) + y_{32}(5A - 3B - 30C)$, $\overline{T}_{y3} = -x_{34}(3A - 5B - 30C) - x_{32}(5A - 3B - 30C)$, $\overline{Z}_4 = -36A + 36B + 180C$, $\overline{T}_{x4} = -y_{34}(3A + 5B + 30C) - y_{41}(5A + 3B + 30C)$, $\overline{T}_{y4} = y_{34}(3A - 5B - 30C) + y_{32}(5A - 3B - 30C)$, $A = (x_{34}y_{21} - x_{21}y_{34})/4$, $B = (x_{32}y_{41} - x_{41}y_{32})/4$, $C = (x_{31}y_{42} - x_{42}y_{31})/4$

## Appendix 2

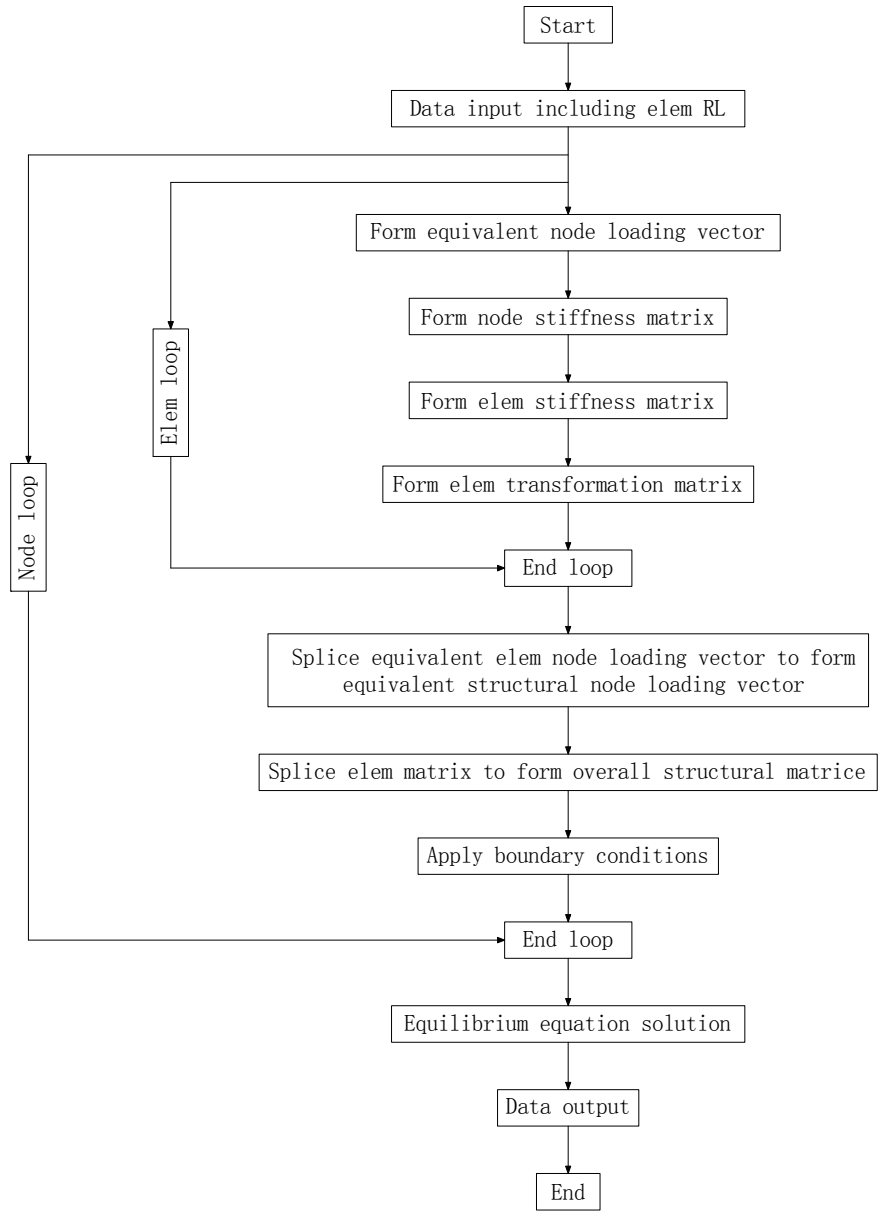

Code flow chart